\documentclass[hidelinks]{article}
\usepackage{times}
\usepackage{amsmath,amsfonts}
\usepackage{graphicx,color}
\usepackage{comment}

\usepackage{listings}
\newcommand{\dealii}{{\textsc{deal.II}}}
\newcommand{\pfrst}{{\normalfont\textsc{p4est}}}

\newcommand{\aspect}{\textsc{ASPECT}}

\newtheorem{remark}{Remark}

\begin{document}

\markboth{Gassm{\"o}ller, Heien, Puckett, and Bangerth}{Particle methods for
  parallel computations}
\title{Flexible and scalable particle-in-cell methods for massively parallel computations}

\author{
  RENE GASSMOELLER\\
  Colorado State University\\  \texttt{rene.gassmoeller@mailbox.org}
  \and
  ERIC HEIEN\\  
  University of California, Davis
  \and
  ELBRIDGE GERRY PUCKETT \\
  University of California, Davis
  \and
  WOLFGANG BANGERTH\\
  Colorado State University
}

\maketitle

\begin{abstract}
  Particle-in-cell methods couple mesh-based methods for the solution
  of continuum mechanics problems, with the ability to advect and
  evolve particles. They have a
  long history and many applications in scientific computing. However,
  they have most often only been implemented for either sequential
  codes, or parallel codes with static meshes that are statically
  partitioned. In contrast, many mesh-based codes today use
  adaptively changing, dynamically partitioned meshes, and can scale
  to thousands or tens of thousands of processors.

  Consequently, there is a need to revisit the data structures and
  algorithms necessary to use particle methods with modern, mesh-based
  methods. Here we review commonly encountered requirements of
  particle-in-cell methods, and describe efficient ways to implement
  them in the context of large-scale parallel finite-element codes that use dynamically changing
  meshes. We also provide practical experience for how to
  address bottlenecks that impede the efficient implementation of these
  algorithms and demonstrate with numerical tests both that our algorithms
  can be implemented with optimal complexity and that they are
  suitable for very large-scale, practical applications. We provide
  a reference implementation in \aspect{}, an open source code for geodynamic
  mantle-convection simulations built on the \dealii{} library. 
\end{abstract}

\section{Introduction}
\label{sec:introduction}

The dominant methodologies to numerically solve fluid flow problems
today are all based on a continuum description of the fluid in the
form of partial differential equations, and
include the finite element, finite volume, and finite difference
methods. On the other hand, it is often desirable to couple these
methods with discrete, ``particle'' approaches in a number of
applications. These include, for example:

\begin{itemize}
\item \textit{Visualization of flows:} Complex, three-dimensional flow
  fields are often difficult to visualize in ways that support
  understanding of their critical aspects. Streamlines are sometimes used
  for this but are typically based on instantaneous flow fields. On
  the other hand, showing the real trajectories of advected particles
  can provide real insight, in much the same way as droplets are used
  for visualization in actual wind tunnel tests.
\item \textit{Tracking interfaces and origins:} In many applications,
  one wants to track where interfaces between different parts of the
  fluid are advected, or where material that originates in a part of
  the domain ends up as a result of the flow. For example, in
  simulating the flow in the Earth's mantle (the application for which
  we developed the methods discussed in this contribution), one often
  wants to know whether the long term motion leads to a complete
  mixing of the mantle, or whether material from the bottom of the
  mantle stays confined to that region.
\item \textit{Tracking history:} In other cases, one wants to track
  the history of a parcel of fluid or material. For example, one may want
  to integrate the amount and direction of strain a rock grain experiences 
  over time as it is
  transported through the Earth's mantle, possibly relaxed over time by
  diffusion processes. This may then be used to understand seismic
  anisotropies, but accumulated strain can also be used in damage
  models to feed back into the viscosity tensor and thereby
  affect the flow field itself.
\end{itemize}

Many of these applications could be achieved by advecting additional,
scalar or tensor-valued fields along with the flow. For example,
interfaces and origins could be tracked by additional ``phase'' fields with
initial values between zero and one, and the history of a volume of fluid could be tracked by additional
fields whose values change as a function
of properties of the flow field or their location. 

On the other hand, advected fields also have numerical disadvantages such as 
undesirable levels of diffusion and failure to preserve appropriate bounds in regions 
with sharp gradients; i.e., overshoot and undershoot.
For these reasons -- and, importantly, simply because some areas of computational
science have ``always'' done it that way -- approaches in which particles
are advected along with the flow are quite popular for some
applications. Use cases and discussions of computational methods can,
for example, be found as far back as 1962 in \cite{Harlow_1962} and are often referred to as
\textit{particle-in-cell} (PIC) methods to indicate that they couple
discrete (Lagrangian) particles and (Eulerian) mesh-based continuum methods.

Many implementations of such methods can be found in the
literature \cite{poliakov1992,gerya2003,McNamara2004,Popov2008,TMK14}. 
However, almost all of these methods are restricted to either simple, structured meshes 
and/or sequential computations. 
This no longer reflects the state of the art of fluid flow
solvers. The current contribution is therefore a comprehensive
assessment of all of the challenges arising when implementing
particle methods in the context of modern computational fluid dynamics (CFD)
solvers. Specifically, we aim at developing methods that also work in
the following two situations:
\begin{itemize}
\item \textit{Unstructured, fully adaptive, dynamically changing
  2d and 3d meshes:} We envision solving complex and dynamically
  changing flow problems in geometrically non-trivial domains. These
  are most efficiently solved on dynamically adapted meshes that use
  mesh refinement and coarsening.
\item \textit{Massively parallel computations:} The methods we aim to
  develop need to scale to computations that run on thousands of cores
  or more, tens of millions of cells or more, and billions of particles.
\end{itemize}
Supporting these two scenarios in mixed particle-mesh methods requires
fundamentally different approaches to the design of data structures
and algorithms. Within this paper, we will therefore investigate all
aspects of implementing such methods required for actual, realistic
applications. Specifically, we will consider the following components,
along with a discussion of their practical performance in typical use cases:
\begin{itemize}
\item Appropriate data structures for the use cases outlined above;
\item Generation of particles;
\item Advection of particles by integrating their trajectory and properties;
\item Treatment of particles as they cross cell and processor
  boundaries in parallel computations;
\item Treatment of particles when the mesh is refined or coarsened adaptively, including
  appropriate load balancing;
\item How to deal with ``properties'' that may be used to characterize
  auxiliary state variables associated with each particle;
\item Methods to interpolate or project particle properties onto the mesh.
\end{itemize}

Our goal in this paper is \textit{to provide a comprehensive
  assessment of all of the steps one needs to address to augment
  state-of-the-art computational fluid dynamics solvers with efficient
  and scalable particle schemes}. While our discussions will be
generic and independent of any concrete software package, we provide a
reference implementation as part of the mantle convection modeling
code \aspect{} \cite{KHB12,aspectmanual}. We will use this reference
implementation to also demonstrate practical properties of our
approaches and thereby validate their usefulness, efficiency,
accuracy, and scalability. In particular, we will show numerical
results up to many thousands of cores, using billions of particles,
thus demonstrating applicability to real-world testcases.

This paper is structured as follows: Section \ref{sec:computational-methods}
presents an overview of the algorithms and data structures we use, including
considerations of algorithmic complexity. Section
\ref{sec:results} then evaluates these methods in scaling
studies and shows that our algorithms fulfill the claimed properties. Section
\ref{sec:application} presents a typical application in geodynamic modeling
that makes use of the novel properties of our algorithm. We conclude in
Section~\ref{sec:conclusion}.

\section{Computational methods}
\label{sec:computational-methods}

In this section, let us discuss the algorithmic challenges listed in the introduction.
Specifically, we address particle generation (Section~\ref{subsec:generation}),
advection (Section~\ref{subsec:advection}), transport between cells and parallel 
subdomains (Section~\ref{subsec:transport}), handling of cells upon mesh refining and 
coarsening (Section~\ref{subsec:adaptive_mesh_handling}), treatment of particle
properties (Section~\ref{subsec:properties}), transfer of information from
particles to mesh cells (Section~\ref{subsec:mapping_properties}), and the
generation of graphical or text-based output of particle data (Section~\ref{subsec:output}). 
In most of these cases, the implementation must allow for different choices of specific 
algorithms; for example, the particle advection algorithm may be the single-step 
forward Euler method or the two or four-step second and fourth-order Runge-Kutta 
methods. 
We will discuss how one can design generic implementations for such cases
and comment on general software design choices for such a framework in
Section~\ref{subsec:implementation}.

The theoretical assessment of performance criteria of algorithms relies
on the use of appropriate data structures. Consequently, we
will start our discussion in Section~\ref{sec:data-structures} with a
presentation of the fundamental storage scheme we propose to use. The
evaluation of algorithmic complexities is then based on the complexity
of accessing elements within these data structures.

\subsection{Data structures}
\label{sec:data-structures}

To evaluate the complexity
of the algorithms in this publication, we assume data structures described in this subsection, which satisfy the following, basic requirements as part of a
mesh-based CFD code:
\begin{itemize}
\item At any given time, each particle is associated with a cell $K$ within
  which it is located at that time.
\item Each cell $K$ contains a number $N_{\text{particles},K}$ particles that may be different
  from cell to cell and from (sub-)time step to (sub-)time step due to
  particles being advected from one cell to another.
\item Each cell $K$ represents a part $\Omega_K$ of the global model domain $\Omega$
\item Each processor $P$ ``owns'' a number $N_{\text{cells},P}$ of cells
  that may change from time step to time step due to mesh refinement and load
  balancing. The number of all particles in these $N_{\text{cells},P}$ cells is
  $N_{\text{particles},P}$.
\item Each particle stores its current location.
\item Each particle stores a global index $i$ that uniquely identifies it.
\item Each particle stores a number $N_{\text{properties}}$ of scalar
  properties.%
\footnote{In practice, these scalar properties may be interpreted as the
    components of a higher-rank tensor such as an accumulated strain. However,
  how this data is \textit{semantically interpreted} is immaterial to our
  discussions here.}
  $N_{\text{properties}}$ is a run-time constant but
  is, in general, not known at compile time.
\item Algorithms must be able to:
  \begin{itemize}
  \item efficiently identify all particles located in a cell, 
  \item efficiently identify the cell for a given particle,
  \item efficiently transfer particles from one cell to another,
  \item efficiently transfer particles from one processor to another,
  \item efficiently evaluate field-based variables at the location of a particle.
  \end{itemize}
\end{itemize}

These requirements make it clear that fixed-sized arrays would be a
poor choice. Rather, we use the following, C++-style structure to represent the
data associated with a single particle:
\begin{lstlisting}[frame=single,basicstyle=\footnotesize]
  template <int dim> struct Particle
  {
    ParticleId   id;
    Point<dim>   location;
    Point<dim>   reference_location;
    double      *properties;
  };
\end{lstlisting}
Here, \texttt{ParticleId} is a type large enough to uniquely index all
particles in a simulation (we provide the choice of a 32-bit or 64-bit unsigned
integer type as a compile time option), and \texttt{Point<dim>} is a data type
that represents the
\texttt{dim} coordinates of a location in ${\mathbb R}^\text{dim}$.
\texttt{location} and \texttt{reference\_location} denote the location of the
particle within the domain $\Omega$ in a global coordinate system, and within
local coordinate system of the cell $K$ that encloses the particle, respectively.
\texttt{properties} points to a dynamically managed memory address
that can store $N_{\text{properties}}$ scalars; this location may be provided
by a memory pool class that manages memory in fixed increments of
$N_{\text{properties}}$ scalars.

Every processor then stores all of its particles in an object declared as
follows:
\begin{lstlisting}[frame=single,basicstyle=\footnotesize]
  std::multimap<cell_iterator, Particle<dim> > particles;
\end{lstlisting}
The multi-map allows the storage of $0$ or more particles for each
cell among
that subset of cells of a mesh that a processor ``owns'', where cells
are represented as iterators into the mesh data structure. 
Concretely, in an adaptive mesh structure that allows for dynamically changing meshes, 
a cell iterator can be represented by a pointer to the mesh object, the
number of refinements a cell has undergone starting at the coarse mesh
(the ``level'' of a cell), and the ``index'' of this cell within the
set of all cells on this level. Since the mesh object will always be
the same, we store the pointer externally, and the iterator is
characterized only by the level and index. These pairs clearly allow
for a lexicographical total ordering and therefore can serve as keys
into a map or multi-map.

\begin{remark}
  \label{rem:ghost-particles}
  For some applications, it is also necessary to store but (usually) neither advect 
  nor otherwise update particles that are located in ``ghost'' cells, i.e.,
  cells that surround the ones owned by one processor, but are owned by 
  another processor.
  We will give examples of this in Section~\ref{subsec:mapping_properties}.
  Whether these particles are stored in the same multi-map, or a separate one, is
  unimportant; however, it is convenient to use the same fundamental data structure.
\end{remark}

We chose a C++ standard template library \texttt{std::multimap}
data structure since it can be
efficiently implemented as a binary search tree, connecting individual
keys (the cell iterators in our case) with multiple values per key
(the particles residing in a particular cell). This container keeps
the particles sorted by their containing cells at all times, and
allows us to efficiently iterate over all cells, handling all
particles in one cell at a time. Each insert, deletion, and search
for individual particles is ${\cal O}(\log(N_{\text{cells},P}))$ in complexity; already presorted 
collections of particles can be merged in ${\cal O}(N_{\text{particles},P})$ complexity, and loops over all cells are of ${\cal O}(N_{\text{cells},P})$ complexity.

The advantages of using a dynamic data structure such as a multi-map indexed by the 
cells are as follows:
\begin{itemize}
\item It is easy to move particles from one cell to another, as well as to remove them 
    from one processor and move them to another.
\item It is easy to loop over all cells, identify all of the particles in the current cell, and loop over them in order to evolve their positions and properties with the velocities and 
    other values defined by the mesh-based flow solver.
\end{itemize}
An alternative is to use an array of particles that individually store pointers to the 
cell that each particle is in.
However, this approach either requires active management of used/unused array 
locations when moving particles between processors, or compressing arrays after every 
particle transfer.
It also either requires linear
searches for particles located in a particular cell, or periodic
sorting of arrays based on the cell iterator key. A careful evaluation of the
operations one typically has to do in actual codes shows that the overhead of
using array-based data structures outweighs their simplicity.

Having thus fixed a data structure, we will be able to describe
concrete algorithms in the following sections and analyze their
complexities based on the complexity of accessing, inserting, or
deleting elements of the multi-map.

\begin{remark}
  \label{rem:reference-coordinates}
  It is often necessary to also know the coordinates $\hat{\mathbf x}$ of the $k$th 
  particle in the \textit{reference coordinate system} of the cell, $K$, in which the 
  particle  currently lies; i.e., to know
  $\hat{\mathbf x}_k = \Phi_K^{-1}(\mathbf x_k)$ where  the (possibly nonlinear) 
  function $\Phi_K : \hat K \mapsto K$ maps the reference cell, $\hat K$,
  (typically, simplices or hypercubes) to a cell $K$ in the physical  domain.
  This can of course be computed on the fly every time it is needed, but since
  evaluating $\Phi_K^{-1}$ is expensive, we also store the
  \textit{reference location} $\hat{\mathbf x}_k$ within the \texttt{Particle} data structure and keep it in sync
  with the global location of the particle when a particle is generated, every time it is moved, or its containing cell is refined or coarsened.
\end{remark}

\subsection{Generation}
\label{subsec:generation}

\begin{figure}
\centering
\includegraphics[width=0.6\textwidth]{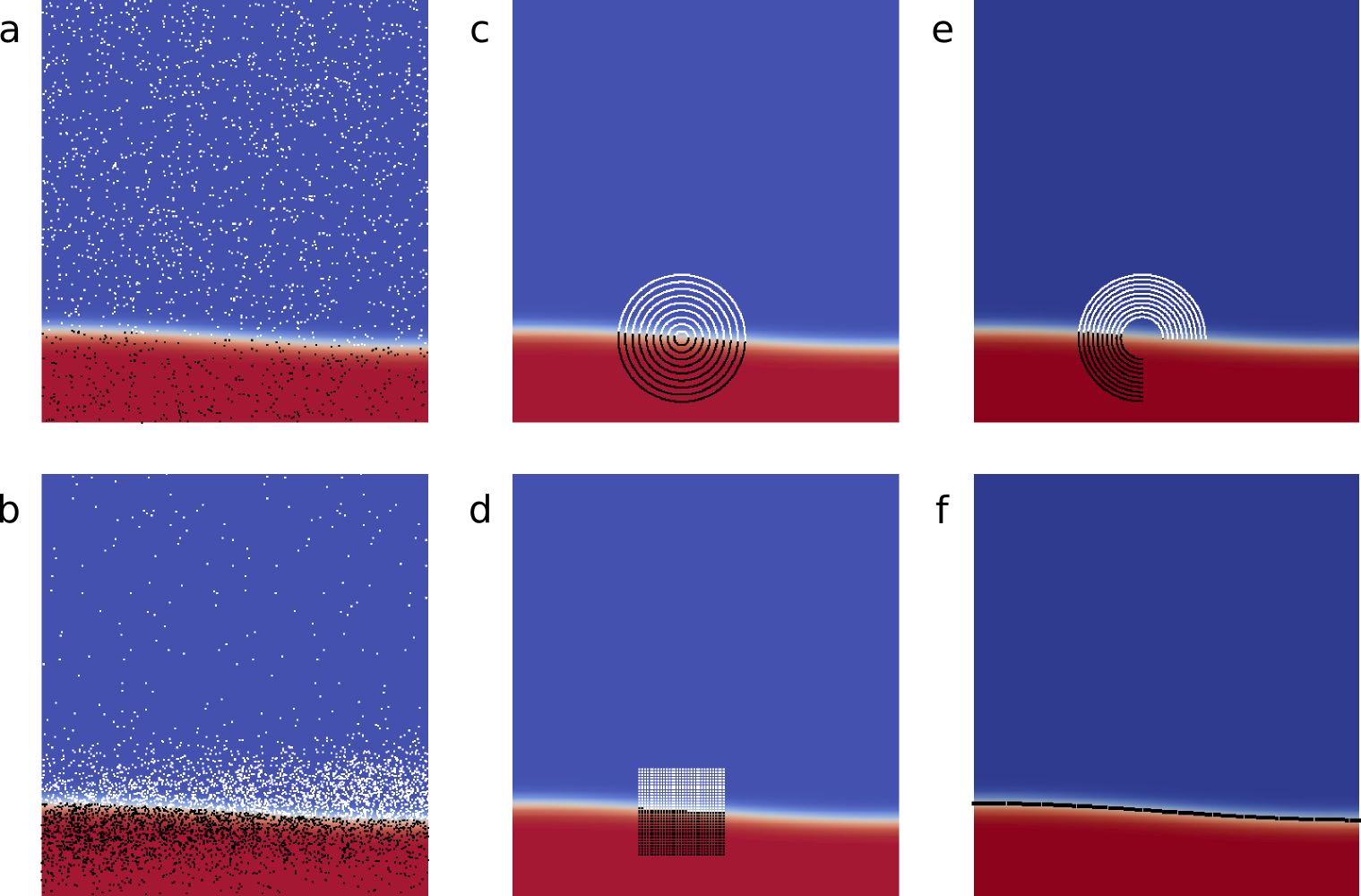}
\caption{\it Common methods to choose the initial particle positions:
    A statistical distribution with (a, b)~uniform or analytically prescribed 
    particle densities, or (c, d)~regular distributions in spherical or box patterns.
    (e, f)~More complicated initial particle patterns can be created
    by variations of the regular distributions or with a user-given
    input file. 
    Model domain and particle coloring is for illustration purposes only and was chosen 
    for a computation that involves the tracking of an interface in which the color 
    indicates the initial position of the particle.}
\label{fig:generation}
\end{figure}

The first step in using particles in mesh-based solvers is their
creation on all involved processors. Depending on their intended uses,
particles may initially be distributed randomly (possibly based on a
probability distribution) or in specific patterns. We will describe
the algorithms necessary for both of these cases in the
following. Examples of use cases are shown in Fig.~\ref{fig:generation}.

\paragraph*{Random particle positions} Randomly chosen particle
locations are often used in cases where particles represent the values
of a field; e.g., the history of material, or origin and movement of a specific type of material.
In these cases, one is often not interested in prescribing the exact location of each 
particle, but randomly chosen locations are acceptable.
The probability distribution, $\rho(\mathbf x)$, from which locations are
drawn would usually be chosen as uniform in the first of the cases
above, and is often the characteristic function of a subdomain in the
second case. Alternatively, one can use a higher particle density in
certain regions, for example to better resolve steep gradients. 
In practice, it is often not possible to ensure that
$\rho(\mathbf x)$ is normalized to $\int_\Omega \rho(\mathbf x) = 1$
where $\Omega$ is the computational domain; thus, the algorithm below
is written in a way so that it does not require a normalized
probability distribution.

To generate $N$ particles on $P$ processors, we use the following
algorithm, running on each processor:
\begin{lstlisting}[frame=single,basicstyle=\footnotesize]
  // compute this processor's weight with regard to %*$\rho(\mathbf x)$*)
  double              local_weight = 0;
  std::vector<double> accumulated_cell_weights(%*$N_\text{cells,P}$*));
  for (cell %*$K$*) in locally owned cells)
    {
      local_weight += %*$\int_K\rho(\mathbf x) \; \text{d}x$*)
      accumulated_cell_weights[%*$K$*)] = local_weight;
    }
  double global_weight = %*$\sum_{p \in \text{processors}} \text{local\_weight}[p]$*);

  // compute this processor's number of particles, 
  // and the id of the first local particle
  ParticleId %*$N_\text{particles,P}$*) = %*$N_\text{particles}$*) * local_weight / global_weight;
  ParticleId local_start_id = %*$\sum_{p=0}^{\text{local\_process\_id}-1} N_\text{particles,p}$*);
  
  // determine in which cells to generate particles
  std::vector<unsigned int> particles_in_cell(%*$N_\text{cells,P}$*));
  for (%*$i \in [0,N_\text{particles,P}$*)))
    {
      random_number = local_weight * uniform_random(0,1);
      CellId %*$K$*) = accumulated_cell_weights.lower_bound(random_number);
      ++particles_in_cell[%*$K$*)];
    }
    
  // create particles in each cell
  ParticleId next_id = local_start_id;
  for (cell %*$K$*) in locally owned cells)
    {
      for (%*$i \in [0,\text{particles\_in\_cell}[K]$*)))
        {
          Particle p;
          p.location           = ...;    // compute a location in K
          p.reference_location = ...;
          p.id                 = next_id;
          particles.insert (pair(%*$K$*), p));
 
          ++next_id;
        }
    }
\end{lstlisting}

Apart from the two global reductions to determine the global weight and
the local start index, all of the operations above are
local to each processor. The overall run time for generating particles
is proportional to $\max_{1\le p\le P} N_{\text{particles},p}$, i.e.,
of optimal complexity in $N$ and, if the number of particles per
processor is well balanced,%
\footnote{However, this is often not the case in practice, see Section~\ref{subsec:adaptive_mesh_handling}.}
 also in $P$. The insertion of each
particle into the multi-map in the algorithm above is, in general, not
possible in constant time. However, since we create particles by
walking over cells in the order they are sorted, it is easy
to modify the algorithm slightly to first put cell-particle pairs
into a linear array (where they will already be sorted by cell) and then
fill the multi-map with all particles at once; this last step can be done with ${\cal
  O}(N_{\text{particles},p})$, and therefore optimal, complexity.

There remains the task of computing a random location in cell $K$,
i.e., the step not yet described in the algorithm above. This step is
not trivial for cells that are not rectangles (in 2d) or boxes (in
3d), in particular if the cell has curved boundaries. In our
implementation, we first determine the axes-parallel bounding box $B_K
\supseteq K$ of $K$, and then keep drawing particles at uniformly
distributed random locations in $B_K$ until we find one that is inside
$K$. This is sufficient for most cases, but has several drawbacks:
(i)~The density of generated particles is uniform in $K$, even if 
$\rho(\mathbf x)$ is not; in particular, even if $\rho(\mathbf x)$ is
zero in parts of the cell (but not everywhere in $K$, so that the cell's
weight is non-zero), then the algorithm may still draw particle
locations in those parts. (ii)~The algorithm becomes inefficient if
the ratio of the volume of $K$ to the volume of $B_K$ becomes small,
because the number of particle locations one needs to draw before
finding one in $K$ is, on average, $|B_K|/|K|$. An example for this situation would
be pencil-like cells oriented along the coordinate system
diagonals. (iii)~Determining whether a trial point $\mathbf x_t$ is
inside $K$ is generally expensive unless $K$ is a (bi-, tri-)linear
simplex because it requires the iterative inversion of the non-linear
mapping $\Phi_K$ from reference cell $\hat K$ to $K$ in order to
determine whether $\Phi_K^{-1}(\mathbf x_t)\in \hat K$.

An algorithm that addresses all of these shortcomings uses a
Metropolis-Hastings (MH) Monte Carlo method to draw a sequence of
locations $\{\hat {\mathbf x}_k\}_{k=1}^{\tilde N_K}$ in the reference
cell $\hat K$ based on the (non-normalized) probability density
$\hat\rho(\hat{\mathbf x}):=\rho(\Phi_K(\hat {\mathbf x}_k))
\; |\det \hat\nabla\Phi_K(\hat {\mathbf x}_k)|$. The proposal
distribution for the MH algorithm can 
conveniently be chosen uniform in $\hat K$, necessitating the second
factor in the definition of $\hat\rho(\hat{\mathbf x})$. This method has the
advantage that it only requires computing the cheap (often polynomial)
forward mapping $\Phi_K$ rather than its inverse, and that it does not
draw points at locations where $\rho(\Phi_K(\hat {\mathbf x}))=0$. On
the other hand, the MH method yields a sequence of points in which
subsequent points may be at the same location. Consequently, to
generate $N_K$ particle locations for cell $K$ that continue to have
the desired probability density, we have to create a longer sequence
of length $\tilde N_K \ge N_K$ and take every $J^\text{th}$ element
where $J$ is larger than the length of the longest subsequence of
equal samples.

We note that because the number of particles per processor is fixed initially based
on $\rho(\mathbf x)$, our algorithm is not \textit{entirely} random in the distribution.
However, in practice we find this does not matter if there are
sufficiently many particles; a completely random
method would then unnecessarily add complexity.

\paragraph*{Prescribed particle locations} An alternative to the
random arrangements of particles are cases where users want to exactly
prescribe initial particle locations. These locations may either be
programmatically described (e.g., on a regular grid within a small
box), or simply be coordinates listed in a file. Maybe surprisingly,
this general case turns out to be more computationally expensive than randomly
generated particle locations.

In either case, let us assume that the initial positions of all
particles are given in an array $\{\mathbf x_k\}$, $k=1\ldots N$. Then
the following algorithm performs particle creation and insertion:
\begin{lstlisting}[frame=single,basicstyle=\footnotesize]
  ParticleId next_id = 0;
  for (%*$k \in [1,N]$*))
    for (cell %*$K$*) in locally owned cells)
      if (%*$\mathbf x_k \in K$*))
        {
          Particle p;
          p.location           = %*$\mathbf x_k$*)
          p.reference_location = %*$\Phi^{-1}_K(\mathbf x_k)$*)
          p.id                 = next_id;
          particles.insert (pair(%*$K$*), p));

          ++next_id;
        }
\end{lstlisting}

In the worst case, the complexity of this algorithm is ${\cal
  O}(N_{\text{particles}}N_{\text{cells},p})$ on processor $p$. This is because, with
general unstructured meshes, we can not predict whether a given
particle's location lies inside the set of cells stored by this
processor without searching through all cells. This limits the
usefulness of the algorithm to moderate numbers of particles. (The
number of cells per processor, $N_{\text{cells},p}$, is often already
limited to at most a few hundred thousands for other reasons.)
However, the algorithm can be accelerated by techniques such as
checking whether a particle's location lies inside the bounding box of
the current processor's set of locally owned cells, before checking
every one of the locally owned cells. 

On hierarchically refined meshes, one can sometimes also find the cell $K$
for a given particle position $\mathbf x_k$ by finding the coarse
level cell in which it is located, and then recursively searching through its
children. This reduces the complexity to ${\cal
  O}(N_{\text{particles}}\;\log N_{\text{cells},p})$. However, it only
works if child cells occupy the same volume as their parent cell; this
condition is often not met when using curved geometries.

\begin{remark}
  \label{rem:local-generation}
  In the paragraph above we assume that the particle positions are known in the
  global coordinate system, and $\Phi_K^{-1}$ has to be evaluated in order to
  find the surrounding cell. If however, the particle coordinates are known in
  the \textit{local} coordinate system of each cell (e.g., if one
  wants to create one particle at each of the cells' midpoints), then
  the above algorithm is much simpler. A loop over all cells and all local
  particle coordinates will generate the particles in the optimal order to
  insert them into the multi-map, and will only involve a (cheap) evaluation
  of $\Phi_K$.
\end{remark}

\subsection{Advection}
\label{subsec:advection}

The key step -- as well as typically the most expensive part -- of any
particle-in-cell method is solving the equation of motion for each particle
position $\mathbf{x_k}=\mathbf{x_k}(t)$,
\begin{equation}
\label{eq:position-ode}
\frac{d}{dt} \mathbf{x}_k (t) = \mathbf u(\mathbf{x},t),
\end{equation}
where $\mathbf{u}(\mathbf{x},t)$ is the velocity field of the
surrounding flow, which is usually computed with a discretized continuum mechanics 
model. 
In practice, the exact velocity $\mathbf u(\mathbf x,t)$ is not available, but only a 
numerical approximation $\mathbf u_h(\mathbf x,t)$ to $\mathbf u(\mathbf x,t)$. 
Furthermore, this approximation is only available at discrete time steps,
$\mathbf u_h^n(\mathbf x) = \mathbf u_h(\mathbf x,t^n)$ and
these need to be interpolated between time steps if the
advection algorithm for integrating~\eqref{eq:position-ode} requires one or more 
evaluations at intermediate times between $t^n$ and $t^{n+1}$.
If we denote this interpolation in time by $\tilde{\mathbf u}_h(\mathbf x,t)$
where $\tilde{\mathbf u}_h(\mathbf x,t^n)= \mathbf u^n_h(\mathbf x)$, then the equation 
the differential equation solver really tries to solve is
\begin{equation}
\label{eq:position-ode-approx}
\frac{d}{dt}\tilde{\mathbf x}_k(t) = 
\tilde{\mathbf u}_h(\mathbf x_k(t),t).
\end{equation}

\begin{remark}
\label{remark:convergence}
Assessing convergence properties of an ODE integrator -- for
example to verify that the RK4 integator converges with
fourth order -- needs to take into account that the error in particle
positions may be dominated by the difference $\mathbf u-\tilde{\mathbf
  u}_h$, instead of the ODE solver error. If, for example, we denote
by $\tilde{\mathbf x}_{k,h}(t)$
the numerical solution of \eqref{eq:position-ode-approx}, then the
error will typically satisfy a relationship like
$$
\| \tilde{\mathbf x}_k(T) - \tilde{\mathbf x}_{k,h}(T) \|
\le
C(T) {\Delta t}_\text{p}^q
$$
where ${\Delta t}_\text{p}$ is the time step used by the ODE solver
(which is often an integer fraction of the time step ${\Delta t}_{\textbf{u}}$ used to
advance the velocity field $\mathbf u$), $q$ the convergence order
of the method, and $C(T)$ is a (generally unknown) constant
that depends on the end time $T$ at which one compares the
solutions. On the other hand, one would typically compute ``exact''
trajectories using the \textit{exact} velocity,
and then assess the error as
$\| \mathbf x_k(T) - \tilde{\mathbf x}_{k,h}(T) \|$.
However, this quantity will, in the best case, only satisfy an
estimate of the form
$$
\| \mathbf x_k(T) - \tilde{\mathbf x}_{k,h}(T) \|
\le
C_1(T) {\Delta t}_\text{p}^q
+ C_2(T) \| \mathbf u-\mathbf u_h \|
+ C_3(T) \| \mathbf u_h-\tilde{\mathbf u}_h \|,
$$
with appropriately chosen norms for the second and third
term. These second and third terms typically converge to
zero at relatively low rates (compared to the order $q$ of 
the ODE integrator) in the mesh size $h$ and the time step size
${\Delta t}_\mathbf{u}$, limiting the overall accuracy of the ODE integrator.
\end{remark}

Given these considerations, and given that ODE integrators require the
expensive step of evaluating the velocity field $\tilde {\mathbf u}_h$
at arbitrary points in time and space, choosing a simpler, less
accurate scheme can significantly reduce the computation time. In our
work, we have implemented the forward Euler, Runge-Kutta 2 and
Runge-Kutta 4 schemes \cite{HW91}. We will briefly discuss them below
and remark that we have found that using higher order or
implicit integrators does not usually yield more accurate solutions.
For simplicity, we will
omit the particle index $k$ from formulas in the remainder of this section.

In the following, for simplicity in exposition, we will assume that
the ODE and PDE time steps $\Delta t_p,\Delta t_{\mathbf u}$ are
equal. We will therefore simply denote them as $\Delta t$. This is
often the case in practice because the velocity field 
is typically computed with a method that requires a
Courant-Friedrichs-Lewy (CFL) number around or smaller than one,
implying that also particles move no more than by one cell diameter
per (PDE) time step. In such cases, even explicit time integrators for
particle trajectories can
be used without leading to instabilities, and all of the methods below
fall in this category. The formulas in the remainder of this section are,
however, obvious to generalize to cases where $\Delta t_p<\Delta
t_{\mathbf u}$. We will also assume in the following that we have
already solved the velocity field up to time $t^{n+1}$ and are now
updating particle locations from $\mathbf x^n$ to $\mathbf
x^{n+1}$. In cases where one wants to solve for particle locations
\textit{before} updating the velocity field, $\tilde{\mathbf u}_h$ can
be extrapolated beyond $t^n$ from previous time steps.

\paragraph*{Forward Euler}
The simplest method often used is the forward Euler scheme,
\begin{align*}
  \mathbf x^{n+1}
  =
  \mathbf x^{n} + {\Delta t}\;\tilde{\mathbf u}_h(t^{n},\mathbf x^n).
\end{align*}
It is only of first order, but cheap to evaluate and often
sufficient for simple cases.

\paragraph*{Runge-Kutta second order (RK2)} Accuracy and stability can be
improved by using a second order Runge-Kutta scheme. The new particle position
is computed as

\begin{align*}
\mathbf{k}_1 &= \frac{{\Delta t}}{2} \tilde{\mathbf u}_h(t^{n}, \mathbf{x}^{n}),
\\
\mathbf{x}^{n+1} &= \mathbf{x}^n + {\Delta t} \; \tilde{\mathbf u}_h\left(t^n+\frac{\Delta t}{2}, \mathbf{x}^{n} + \frac{\mathbf{k}_1}{2}\right).
\end{align*}

\paragraph*{Runge-Kutta fourth order (RK4)} A further improvement of particle
advection can be achieved by a fourth order Runge-Kutta scheme that computes
the new position as 

\begin{align*}
\mathbf{k}_1 &= {\Delta t} \; \tilde{\mathbf u}_h\left(t^{n}, \mathbf{x}^{n}\right), \\
\mathbf{k}_2 &= \frac{{\Delta t}}{2} \; \tilde{\mathbf u}_h\left(t^n+\frac{\Delta t}{2}, \mathbf{x}^n + \frac{\mathbf{k}_{1}}{2}\right), \\
\mathbf{k}_3 &= \frac{{\Delta t}}{2} \; \tilde{\mathbf u}_h\left(t^n+\frac{\Delta t}{2}, \mathbf{x}^n + \frac{\mathbf{k}_{2}}{2}\right), \\
\mathbf{k}_4 &= {\Delta t} \; \tilde{\mathbf u}_h\left(t^{n+1}, \mathbf{x}^n +
\mathbf{k}_{3}\right), \\
\mathbf{x}^{n+1} &= \mathbf{x}^n + \frac 16 \mathbf{k}_1 + \frac 13 \mathbf{k}_2 + \frac 13 \mathbf{k}_3 + \frac 16 \mathbf{k}_4.
\end{align*}

The primary expense in all of these methods is the evaluation of the
velocity field ${\mathbf u}_h^n$ and ${\mathbf u}_h^{n+1}$ at
arbitrary positions $\mathbf x$. By sorting particles into the cells
they are in, at least we know which cell $K$ an evaluation position
$\mathbf x$ lies in. However, given that the velocity fields $\mathbf
u_h$ we consider are typically finite element fields defined based on
shape functions whose values are determined by mapping a reference
cell $\hat K$ to each cell $K$ using a transformation
$\mathbf x=\Phi_K(\hat{\mathbf x})$, evaluation at arbitrary points requires the
expensive inversion of $\Phi_K$. This step can not be avoided, but by storing the
resulting reference coordinates $\hat{\mathbf x}$ of each point after updating
$\mathbf x$, as discussed in
Remark~\ref{rem:reference-coordinates}, we at least do not have to repeat the
step for many of the algorithms below.

\subsection{Transport between cells and subdomains}
\label{subsec:transport}

At the end of each (particle) time step -- or in the case of the RK2 and RK4
methods, at the end of each substep -- we may find that the updated
position for a particle that started in cell $K$ may no longer be in
$K$. Consequently, for each particle, every (sub-)step is followed by
a determination of whether the particle is still inside $K$. This can be done
as part of computing the reference coordinates $\hat{\mathbf x}_k$ and testing
whether $\hat{\mathbf x}_k \in \hat K$. If the particle has moved out of $K$, we
need to find the new cell $K'$ in which the particle now
resides.

In parallel computations, particles may also cross from a cell owned
by one processor to a cell owned by another processor during an
advection step, or even during an advection substep.
Consequently, ownership of particles needs to be transferred efficiently. 
To avoid the latency of transferring individual particles, it is most
efficient to collect all data 
that needs to be shipped to each destination into a single buffer. In
practice, it is sufficient to implement
communication patterns that cover the exchange of particles between
processes whose locally owned cells neighbor each other, i.e., using
point-to-point messages. This is possible, in particular, if the (ODE)
time step is chosen so that the CFL number is less than
one, because then particles travel no more than one cell diameter in
each step. With few exceptions (see below), particles will
therefore either remain within a cell, or end up in a neighboring cell
that is also either owned by the current process, or a ghost cell whose
owner is known by the current process.

Our reference implementation employs the following algorithm to
update the cell-particle map, executed for
each particle stored on the current processor:
\begin{lstlisting}[frame=single,basicstyle=\footnotesize]
if (particle p has left its surrounding cell %*$K$*))
  {
    %*$K'$*) = find_current_cell_for_particle(p,%*$K$*));              // search-cell
    if (%*$K'\;$*) is locally owned cell)                          // option 1
      remove (p,%*$K$*)) from multi-map, add (p,%*$K'$*)) to multi-map;
    else if (%*$K'$*) is ghost cell on current processor)        // option 2
    {
      transmissions[owner of %*$K'$*)].push_back (p,%*$K'$*));
      remove (p,%*$K$*)) from multi-map;
    }
    else                                                   // option 3
      remove (p,%*$K$*)) from multi-map;
  }

for (neighbor n in neighbors)
  communicate_particles(n,transmissions[n]);
\end{lstlisting}

The vast majority of particles remain in the current cell,
end up in a locally owned cell (option 1), or a cell
owned by another process (option 2). In almost all cases where a
particle moves out of its current cell $K$, the cell it
ends up in is in fact a neighbor of $K$.
A small fraction of cases, however, do not fall in these
categories. One example is if the neighbor cell of $K$ is refined, and transport
by one diameter of $K$ results into transport into one of the children
of the neighboring (coarse) cell that is not adjacent to $K$. If this
child is owned by the same process, it can be found and the particle can be
assigned to it (option 1). Otherwise, it is not a ghost cell (because
typically only immediate neighbors of locally owned cells are ghost cells), and
we can not determine its owner. Such particles are then discarded
(option 3). Another example is a particle that
is close to the boundary and that is transported outside the domain
by the ODE integrator; such a particle also needs to be
discarded (again option 3). We have found that even over the
course of long simulations, only a negligible fraction of particles (less than 
0.002\% of all particles per 1,000 RK2 time steps with a CFL number of one)
is lost because of these two mechanisms. As expected, an
explicit Euler scheme increases the number of lost particles significantly 
(1.5\% per 1,000 time steps with identical time step size).
A reduction of the time step size to 0.5 times CFL entirely
avoids all particle losses for RK2, and reduces it significantly for the
forward Euler integrator (0.375\% per thousand time steps).
Avoiding particle loss could be accomplished by further decreasing or
adaptively changing the time step length.  However, given the small
overall loss and added computational expense of a robust solution,
dropping particles that fall out of bounds is a reasonable approach.

The algorithm above requires finding the cell an advected particle
is in now (if any), marked by the \texttt{search-cell} comment. On unstructured
meshes, this in general requires ${\cal 
  O}(N_{\text{cells},P})$ operations; furthermore, because many of the ${\cal
  O}(N_{\text{particles},P})$ cross to a different cell, this step is not of
optimal complexity. On the other hand, in the vast majority of cases, particles
only cross from one cell to its (vertex) neighbors. Consequently, we first
search all of these neighbors, at a cost of ${\cal O}(1)$ per moved
particle since the number of neighbors is typically bounded by a relatively
small constant. Only the very small fraction that do not end up on a
neighbor then requires a complete search over all cells. 

Testing whether a particle is inside a cell $K'$ costs ${\cal
  O}(1)$ operations, but it is expensive. We can accelerate the search
by searching the vertex neighbors of $K$ in an order that makes it
likely that we find the cell surrounding the particle early.
Following some experimentation, we found that the following strategy
works best: Let $\mathbf x$ be the particle's current position,
$\mathbf v$ be the vertex of $K$ closest to $\mathbf x$,
and $\mathbf c_{K'}$ be the center of cell $K'$. Let $\mathbf
a=\mathbf x-\mathbf v$ be the vector from the closest vertex of $K$ to
the particle, and $\mathbf b_{K'}=\mathbf c_{K'}-\mathbf v$ be the
vector from the closest vertex to the center of cell $K'$. Then we
sort the vertex neighbors $K'$ of $K$ by descending scalar
product $\mathbf a \cdot \mathbf b_{K'}$ and search them in this order
whether they contain the particle's new location $\mathbf x$. In other
words, cells with a center in the direction of the
particle movement are checked first. In practice,
we find the new cell in the first try for most cases. On the other
hand, several simpler
criteria -- like the distance between particle and cell center -- fail
more often, in particular for adaptively refined neighbors. 

The deletion and re-insertion of particles into the local multi-map
can be optimized in a similar way to the generation algorithm by first collecting all moved particles in a linear array, and
then inserting them in one loop. This is advantageous because particles of the same
cell tend to move into the same neighbor cells.

The communication of particles to neighboring processes is handled as a two-step process. First, two
integers are exchanged between every neighbor and the current process,
representing the number of particles that will be sent and received
between the respective processes. In a second step every process
transmits the serialized particle data to its neighbors and receives
its respective data from its neighbors. This allows us to implement all
communications as non-blocking point-to-point MPI transfers, only
generating $O(1)$ transmissions and $O(N_{\text{particles},P})$ data per
process. Since we already determined which
cell contains this particle on the old process, we also transmit this information to the new
process, avoiding a repeat search for the enclosing cell.

\subsection{Dealing with adaptively refined, dynamically changing meshes}
\label{subsec:adaptive_mesh_handling}

Over the past two decades, adaptive finite element methods have
demonstrated that they are vastly more accurate than computations on
uniformly refined meshes \cite{Car97,AO00,BR03}. In the current context, such
adaptively refined, dynamically changing meshes present two particular
algorithmic challenges discussed in the following.

\paragraph*{Mesh refinement and repartitioning}
In parallel mesh-based methods, refinement and coarsening typically
happens in two steps: First, cells are refined or coarsened
separately on each process. Particles will then have to
be distributed to the children of their previous cell
(upon refinement), or be merged to the parent of their previous cell (upon coarsening). The second step, after local mesh
refinement and coarsening, requires that the new mesh is
redistributed among the available processes to achieve an efficient
parallel load distribution~\cite{p4est,BBHK10}. During this step, particles need to be
redistributed along with their mesh cells. To keep this process as
simple as possible we append the serialized particle data to other
data already attached to a cell (such as the values of degrees of
freedom located in a cell that correspond to solution vectors, or
vertex locations),  and transmit all data at the same time. We can
therefore utilize the same process typically employed in pure field
based methods, and that is well
supported by existing software for parallel mesh handling \cite{p4est}. In
particular, this approach can use existing and well-optimized
bulk communication patterns, and avoids sending particles individually
or having to re-join particles with their cells.%
\footnote{In the particular implementation upon which we base this
  step in our numerical examples,
  provided by the \pfrst{} library, the data attached to cells to be
  transferred has to have a fixed size, equal for all cells. This is
  not a restriction for mesh based methods that use the same
  polynomial degree for finite element spaces on all cells. However,
  it leads to inefficiencies if the number of particles per cell
  varies widely, as it often does in particle-based simulations. This,
  however, is only a drawback of the particular implementation and can
  easily be addressed in different implementations. As we will show in
  Section~\ref{sec:results}, the inefficiency has no practical impact
  on the overall performance of our simulations.}

\paragraph*{Load balancing}
The mesh repartitioning discussed in the previous paragraph is
designed to redistribute work equally among all
available processes. For mesh-based methods, this typically means
equilibrating the number of cells each process ``owns'', as the
workload of a process is generally proportional to the number of cells
in all important stages of mesh-based algorithms (e.g., assembly,
linear solves, and postprocessing). Consequently, equilibrating the
number of cells between processes also leads to efficient parallel
codes. 

On the other hand, in the context of mesh-particle hybrid methods such
as the ones we care about in this paper, the number of particles per
cell is typically not constant and frequently ranges from zero to a
few dozen or a few hundred. Consequently, rebalancing the mesh so that
each process owns approximately the same number of cells including all
particles located in them, leads to rather unbalanced workloads during
all particle-related parts of the code. Conversely, rebalancing the
mesh so that each process owns approximately the same number of
particles leaves the mesh-based parts of the code with unbalanced
workloads. Both reduce the parallel efficiency of the overall
code.

The only approach to restore perfect scalability is to partition cells
differently for the mesh-based and particle-based parts of the
code. This is possible because one typically first computes the mesh
based velocity field, and only then updates particle locations and
properties. The mesh partitioning step between the two phases of the
overall algorithm then simply follows the outline discussed above. On
the other hand, one can not avoid transporting both mesh and particle
data during these rebalancing steps%
\footnote{Indeed, in many cases this will require transportation of
  not only some, but essentially \textit{all} cell-based and/or particle-based
  data between processors. See for example \cite{p4est}.}
because each phase of the
algorithm requires all data from the other (for example, the particles
may encode material information necessary in computing the velocity
field, whereas the velocity field is necessary to update particle
locations). Consequently, the amount of data that has to be
transported twice per time step is significant.

In practice, some level of imbalance can sometimes be tolerated. In
those cases, one can work with the following compromise solutions:

\begin{itemize}
\item \textit{Repartition mesh to combined weight of particles and
  cells.}
Instead of estimating the workload of each cell during the
rebalancing step as either a constant (as in pure mesh-based methods) or
proportional to the number of particles in a cell (as in pure
particle-based methods), one can estimate it as an appropriately weighted sum of
the two. The resulting mesh is optimal for neither of the two phases
of each time step, but is overall better balanced than either of the
extremes.

\item \textit{Ignore imbalance.} As long as the number of particles is
  small -- where the particle component of the code consequently
  requires only a small fraction of the overall runtime -- one may
  simply ignore the imbalance. A typical case for this is if particles
  are only used to output information for specific points of interest,
  e.g., accumulated strains over the course of a simulation, or to
  track where material from a small set of locations is transported.

\item \textit{Adjust particle density to mesh by particle generation.} 
If the region of interest and highest mesh resolution of a model is
known in advance, one can plan the particle generation accordingly and
align the particle density to the expected mesh density. This is most
useful in cases where, for example, pre-existing interfaces should be tracked
that are expected to only advect but not diffuse over time. This
alignment then not only increases the particle resolution in
regions of interest, but also automatically improves parallel
efficiency and scaling. 

\item \textit{Adjust particle density to mesh by particle population management.}
In cases where the regions of high mesh density are not known in
advance, or in case where particles tend to cluster or diverge in
certain regions of the model, it can be necessary to manage the
particle density actively during the model run. This would include
removing particles from regions with high particle density or adding
particles in regions of low density. If done appropriately, the result
will be a mesh where the average number of particles per cell is managed so
that it remains approximately constant.

\item \textit{Adjust mesh to particle density.}
Instead of adjusting the particle density to align with the mesh
density, the mesh density can also be aligned to the particle
density. This is useful if the feature of interest in a model is most
clearly defined by the particles, e.g., by a higher particle density
close to an interface. As in the previous alternative, the alignment
of mesh and particle density yields better parallel efficiency and scaling.
\end{itemize}

While the last two approaches lead to better scalability, they may of course not
suit the problem one originally wanted to solve. On the other hand, generating
additional particles upon refinement of a cell, and thinning out particles upon
coarsening, is a common strategy in existing codes
\cite{Popov2008,Leng2011}.

\subsection{Properties}
\label{subsec:properties}

Apart from their position and ID, particles are often utilized to
carry various properties of interest. In practice, we have
encountered situations where properties are scalars, vectors, or
tensors; their values may never be updated, may be updated using
finite element field values in every time step, or only when their
values are used to generate graphical output; and they may be
initialized in a variety of ways. Examples for these cases are:
(i) particles initialized by a prescribed function that are used to
indicate different materials in the computational domain; (ii)
particles that store their initial
position to track the movement of material; (iii) particles that represent some part
of the current solution (e.g., temperature, velocity, strain rate) at
their current position for further postprocessing; (iv) particles
whose properties represent the evolution of quantities such as the
accumulated strain, or of chemical compositions.

An example of this last case is a damage model in which a damage variable
increases as material is strained, but also heals over time. This can
be described by the equation
\begin{align*}
  \frac{d}{dt} d_k(t) = \alpha \lVert \dot{\varepsilon}(\mathbf u(\mathbf x_k(t)))\rVert
                       -\beta d_k(t),
\end{align*}
where $\dot{\varepsilon}(\mathbf u(\mathbf x_k(t))) =
       \frac 12 \left(\nabla \mathbf u(\mathbf x_k(t)) +
                      \nabla \mathbf u(\mathbf x_k(t))^T \right)$ is
the strain rate at the location of the $k$th particle, and
$\alpha,\beta$ are material parameters. Similarly, applications in the
geosciences often require integrating the total (tensor-valued) deformation
$F_k$ that particle $k$ has
undergone over the course of a simulation \cite{Hall2000,Becker2003}, leading
to the system of differential equations \cite{McKenzie1983}
\begin{align*}
  \frac{d}{dt} \mathbf{F}_k(t)
  = 
  \nabla \mathbf u(\mathbf x_k(t))
  \mathbf{F}_k(t),
\end{align*}
where $\nabla \mathbf u(\mathbf x_k(t))$ is the velocity gradient tensor,
and $\mathbf{F}_k(0) = \mathbf{I}$.
Clearly, these differential
equations may also reference other field variables than the velocity. In all cases, the
primary computational challenge is the evaluation of field variables
at arbitrary particle positions. This is no different than what
was required in advecting particle positions in
Section~\ref{subsec:advection}, and greatly accelerated by storing the reference
location of each particle in the coordinate system of the surrounding cell
(see Remark~\ref{rem:reference-coordinates}).

In many areas where one uses simulations to qualitatively try to
\textit{understand behavior}, rather than determine a single quantitative
answer, one often wants to track several properties at the same time;
consequently, implementations need to allow assigning arbitrary combinations of
properties -- including their initialization and update features -- to particles.

We have found it useful to only have a single kind of particle in each
simulation, i.e., all particles store the same kinds of properties.
This allows storing information about the semantics of these
properties only once by a property manager object. The manager also
allows querying properties by name, and finding the position of
a particular property within the vector of numbers each particle
stores for all properties. Because the number of scalars that
collectively make up the properties of each particle is a run-time constant, one
can simplify memory management by utilizing a memory pool consisting of a consecutive array
that is sliced into fixed size chunks and that is dynamically managed.

As discussed in the previous section, some models require dynamically
generating or destroying particles during a model run. Particle
properties therefore also need to describe how they
should be initialized in these cases. In practice, we
have encountered situations where one initializes
new particles' properties in the same way as at the start of the computation,
and other cases where it is appropriate to interpolate properties from
surrounding, existing particles.

\subsection{Transferring particle properties to field based formulations}
\label{subsec:mapping_properties}

The previous section was concerned with updating particle properties 
given certain field-based quantities. 
This section is concerned with the opposite problem; namely, how to use properties 
stored on the particles to affect field-based variables.
Examples include damage models, where a damage variable $d_k$ such as the one described 
above would influence the properties of the solid or fluid under consideration.
For example, if the material is a metal, damage increases the material's strength, 
whereas in the case of geological applications, damage typically decreases it.

In field-based FEM methods, one typically evaluates material
properties at the quadrature points. However, particles are not usually
located at the quadrature points and consequently one needs a method
to interpolate or `project' properties from the particle locations to
the quadrature points.
Many such interpolation algorithms, of varying purpose, accuracy, and computational 
cost, have been proposed in the literature~\cite{gerya2003,deubelbeiss2008,TMK14} and 
one's choice of algorithm will typically depend on the intended application.
Therefore, generic implementations of this step need to provide an interface for
implementing different interpolation algorithms.
Interpolation algorithms found in the computational geosciences alone include:
(i)~nearest-neighbor interpolation among all particles in the current cell;
(ii)~arithmetic, geometric, or harmonic averaging of all particles located in the same cell;
(iii)~distance-weighted averaging with all particles in the current cell or within a 
finite distance; 
(iv)~weighted averaging using the quadrature point's shape function values at  the 
particle position as weights;
(v)~a (linear) least-squares projection onto a finite dimensional space.

While these methods differ in their computational cost, most can be implemented with 
optimal complexity since our data structures store all particles sorted by cell, and 
since the algorithms do not require data from other processors.
The exception is the averaging scheme (iii) if the search radius extends past 
the immediate neighbor cells, and if schemes (ii)~and~(iv) are
extended to include particles from such `secondary' neighbor cells.
On the other hand, as long as only information from \textit{immediate} neighboring cells is 
required, these algorithms can be implemented without loss of optimal complexity if the 
procedure discussed in Section~\ref{subsec:transport} is extended to exchange particles
located in one layer of ghost cells for each process.
(See also Remark~\ref{rem:ghost-particles}). 
Methods that require information from cells further away than one layer around a 
processor's own cells pose a significant challenge for massively parallel computations; 
we will not discuss this case further.

\subsection{Large-scale parallel output}
\label{subsec:output}

The last remaining technical challenge for the presented algorithms is
generating output from particle information for postprocessing or
visualization. The difficulty in 
these cases is associated with the sheer amount of data when dealing
with hundreds of millions or billions of particles.

In some cases, particles fill the whole domain with a
high particle density. Outputting \textit{all} of them may then be
unnecessary or even impractical because of the size of the resulting
output files. In these cases, we have found
that it is usually sufficient to instead output an interpolation of
the particle properties onto the mesh, as discussed in the previous
section. 

On the other hand, in applications in which the particles act as
tracers, or in which they store
history-dependent properties with strong spatial gradients, particle
output can provide valuable additional information. In this case, outputting
data for all particles can yield challenging amounts of data (e.g. a single 
output file for one billion particles only carrying a single scalar property 
comprises 36~GB of data in the compressed HDF5 format).
Fortunately, file formats such as 
parallel VTU \cite{vtk} or HDF5 \cite{hdf5} allow each process to write their
own data independently of that of other processes, and in
parallel. Since in practice writing tens of thousands of files per
particle output can be as challenging for parallel file systems as
writing a single very large file, it is often beneficial to group
particle output into a reasonable number of files per time step (say,
10-100 files) or use techniques like MPI~Input/Output to reduce the
number of files created. All steps of this process can
therefore again be achieved in optimal complexity.

\subsection{Implementation choices}
\label{subsec:implementation}

As outlined above, many of the pieces of an overall particle-based
algorithm can be chosen independently. For example, one may want to
use a random particle generation method, and advect them along with a
Runge-Kutta 2 integrator.

\begin{figure}
\centering
\includegraphics[width=0.85\textwidth]{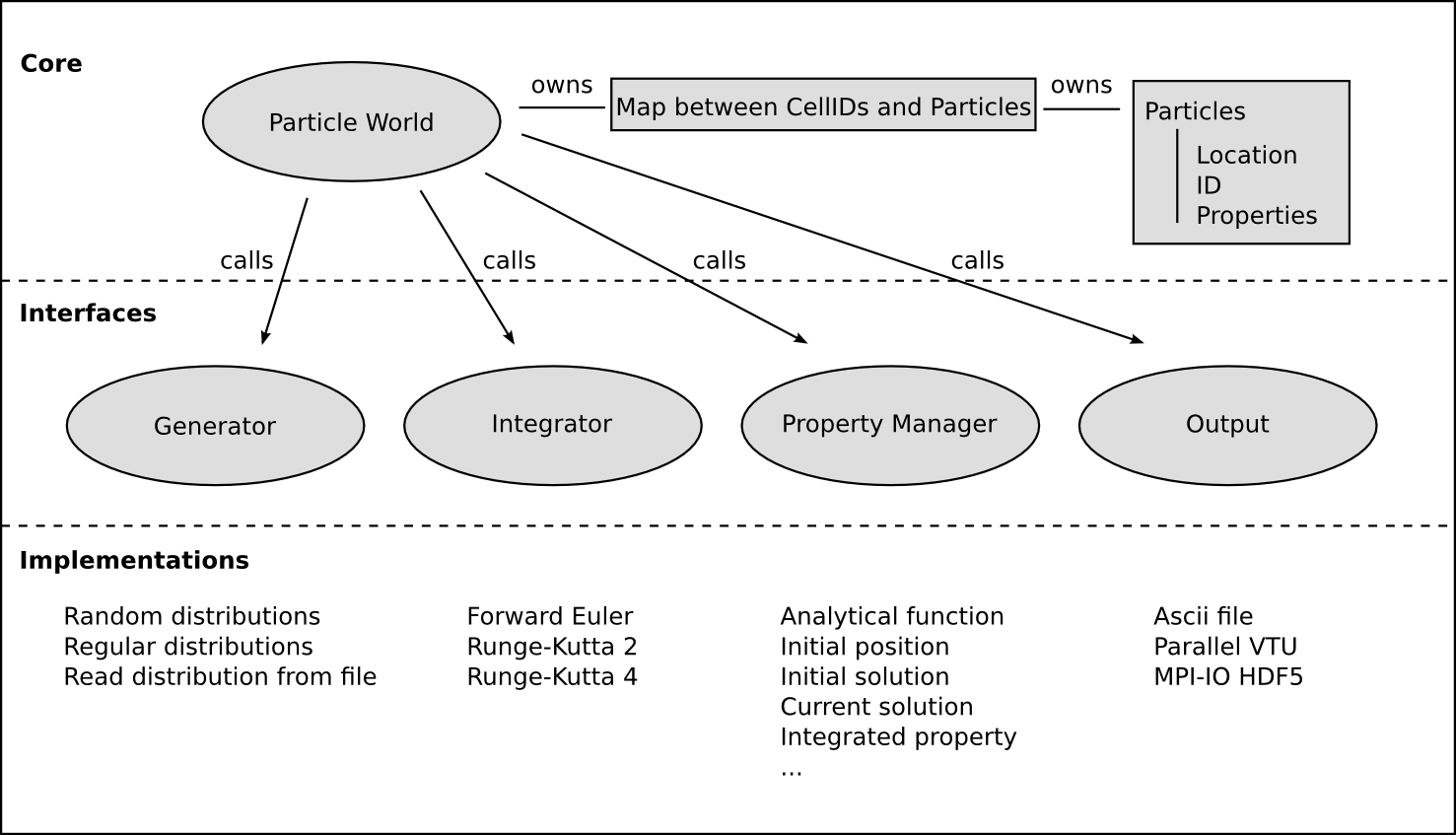}
\caption{\it Structure of the overall class hierarchies used in our reference implementation of the algorithms described in this paper. The functionality is organized in layers which encapsulate the core functionality from specific algorithmic choices. This structure allows flexibly adjusting the algorithm to individual problems at runtime.}
\label{fig:structure}
\end{figure}

This flexibility is conveniently mapped onto implementations in
object-oriented programming languages by defining interface classes for each
step. Concrete algorithms can then be implemented in the form of derived
classes. Fig.~\ref{fig:structure} shows an example of this structure, as
realized in our reference implementation in the \aspect{} code.

\section{Results}
\label{sec:results}

We have implemented the methods discussed in
the previous sections in one of our codes -- the \aspect{} code to simulate
convection in the Earth's mantle \cite{KHB12,aspectmanual}, and that is based on
the deal.II finite element library \cite{BHK07,dealII84} -- with the
goal to provide it in a generic way usable in a wide variety of
settings. To verify our claims of performance and scalability, we
report results of numerical experiments in this section and show that,
as designed, all steps of our algorithms scale well even to very large
problem sizes.%
\footnote{All scaling and efficiency tests were performed on the Cray
  XC-40 cluster at the North-German Supercomputing Alliance (HLRN),
  using nodes with two Xeon E5-2860v3, 12-core, 2.5 GHz CPUs, and a
  Cray Aries interconnect. Each configuration was run 3 times and
  timings were then averaged. We ran each
  configuration for 10 time steps to average over individual time
  steps. When timing individual sections of a program, we introduce
  barriers to ensure accurate attribution to specific
  parts of the algorithm.} 
We also
assessed the correctness of the implemented advection schemes through convergence
tests in spatially and temporally variable flow.

\subsection{Scalability for uniform meshes}
\label{subsec:scalability-uniform}

\begin{figure}[tb]
\centering
\includegraphics[width=.6\textwidth]{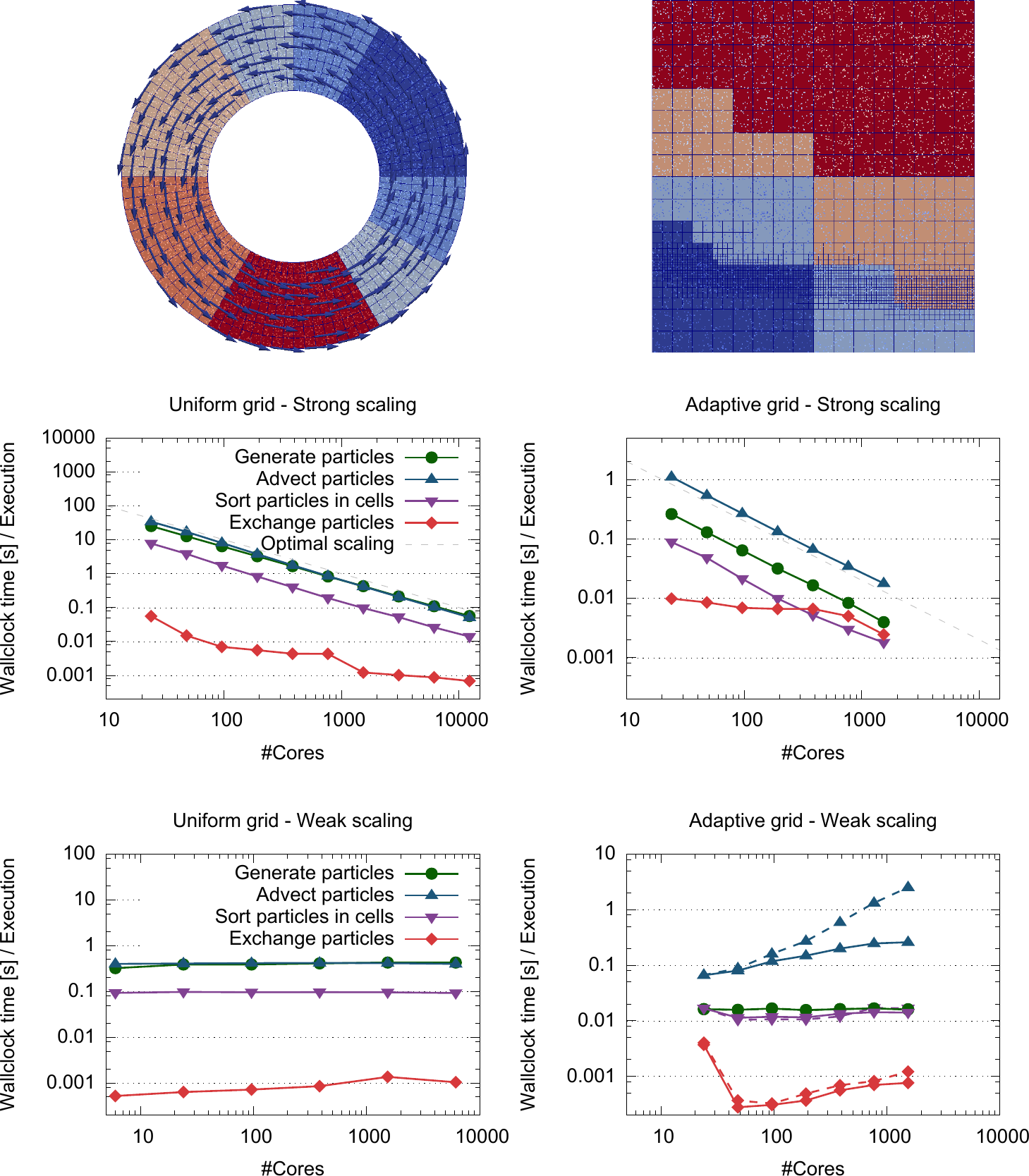}
\caption{\it Scaling of algorithms. Left column: Results for a
  uniformly refined mesh. Right column: Results for an adaptively
  refined mesh. Top row: Model geometry and partition between
  processes. Center row: Strong scaling for a constant number of
  cells and particles and variable number of processes. Bottom row:
  Weak scaling where the numbers of cells and
  particles per process is kept constant. The bottom right panel contains data without load balancing techniques (dashed lines), and with balanced repartitioning as described in Subsection~\ref{subsec:adaptive_mesh_handling} (solid lines). For more information see the main text.}
\label{fig:scaling}
\end{figure}

We first show scalability of our algorithms using a simple
two-dimensional benchmark case with a static and uniformly refined mesh. We
employ a circular-flow setup in a spherical shell, with
no flow across the boundary. Particles
are distributed randomly with uniform density (see
Fig.~\ref{fig:scaling}, top left). 
Because we choose a constant angular velocity of the material, this model creates a much higher
fraction of particles that move into new cells per time step, compared to
realistic applications where large parts of the domain move more slowly
than the fastest regions. Thus, realistic applications
spend less time in the ``Sort particles'' section of the code
(see Section~\ref{subsec:transport}), as we will show in the next testcase
below.

The left column of Fig.~\ref{fig:scaling} shows excellent weak and
strong scaling over at least three orders of magnitude of model
size. For a fixed problem size (strong scaling), we use a mesh with
$1.024\cdot 10^7$ degrees of freedom 
and $1.536\cdot 10^7$ particles. Increasing the number of processes
from 12 to 12,288 shows an almost perfect decrease in wall time for
all operations, despite the rather small problem each process has to
deal with for large numbers of processes. 

Keeping the number of degrees of freedom and particles per core fixed
and increasing the problem size and number of processes accordingly
(weak scaling, bottom left of the figure), the wallclock time stays
constant between 6 and 6,144 processes. In this test each process owns
$6.7\cdot 10^3$ degrees of freedom and $1.0\cdot 10^4$ particles.%
\footnote{Each refinement step leads to four times as many cells, and
  consequently processes. 6,144 cores was the last multiple to which
  we had access for timing purposes.}
Results again show excellent scalability, even to very large problem sizes.

\subsection{Scalability for adaptively refined meshes}
\label{subsec:scalability-adaptive}

Discussing scalability for \textit{adaptive} meshes is more complicated
because refinement does not lead to a predictable increase in the
number of degrees of freedom. We use a setup based on the benchmarks
presented in \cite{JGRB:JGRB11097}. Specifically,
we use a rectangular domain $[0,0.9142]\times[0,1]$ that contains a
sharp non-horizontal interface separating a less dense lower layer
from a denser upper layer. The shape of the interface then leads to a
Rayleigh-Taylor instability. For the strong scaling tests, we adaptively coarsen
a $2,048\times 2,048$ starting mesh, retaining fine cells only in the
vicinity of the interface. The resulting mesh has about
$1.4\cdot 10^6$ degrees of freedom. We use
approximately 15 million, uniformly distributed particles. This setup
is then run on different numbers of processors.

The results in Fig.~\ref{fig:scaling} show that
\textit{strong} scaling for the 
adaptive grid case is as good as for the uniform grid case, nearly decreasing the
model runtime linearly from 12 to 1,536 cores. The model is too
small to allow a further increase in parallelism.

Setting up weak scaling tests is more complicated. Since we can not
predict the number of degrees of freedom for a given number of mesh adaptation steps,
we first run a series of models that all start with a
$256\times256$ mesh and  adaptively refine it a variable
number of times. For each model we assess the
resulting number of degrees of freedom. We then select a series of
these models whose size is approximately a power of two larger than
the coarsest one, and run the
models on as many processes as necessary to keep the ratio between
number of degrees of freedom to number of processes
approximately constant. In practice this number varies between 
41,000 and 48,000 degrees of freedom per process; we have
verified this variability does not significantly affect
scaling results. For simplicity we keep the mesh fixed after the first
time step. Each of the models uses as many particles as degrees of
freedom, equally distributed across the domain.

The weak scaling results are more difficult to interpret than the
strong scaling case. In our
partitioning strategy, we only strive to balance the number of cells
per process. However, because the particle density is constant while
cell sizes strongly vary, the imbalance in the number of particles per
process grows with the size of the model. This is easily seen in the
top right panel of Fig.~\ref{fig:scaling} in which all four processors
own the same number of cells, but vastly different areas and
consequently numbers of particles. Consequently, runtimes for
some parts of the algorithm -- in particular for particle advection --
grow with global model size (dashed lines in the bottom right panel of
Fig.~\ref{fig:scaling}).

As discussed in Section \ref{subsec:adaptive_mesh_handling}, this
effect can be addressed by weighting cell and particle numbers in load
balancing.  The solid lines in the bottom right panel of
Fig.~\ref{fig:scaling} show that with appropriately chosen weights,
the increase in runtime can be reduced from a factor of 60 to a factor
of 3. To achieve this, we introduce a cost factor $W$ for each
particle, relative to the cost of one cell. The total cost of each
cell in load balancing is then one (the cost of the field-based
methods per cell) plus $W$ times the number of particles in this
cell. $W=0$ implies that we only consider the number of cells 
for load balancing, whereas $W=\infty$ only considers the number of particles. 
In practice, one will typically choose $0\le W<1$; the
optimal value depends on the cost of updating
particle properties, the chosen particle advection scheme, and the
time spent in the finite-element solver. For realistic applications,
we found $W=0.01$ to be adequate. On the other hand,
computational experiments suggest that it is not important to
\textit{exactly} determine the optimal value since the overall runtime
varies only weakly in the vicinity of the minimum.

\section{Applications}
\label{sec:application}

We illustrate the applicability of our algorithms using two examples
from modeling convection in the Earth's mantle. The first is based on a
benchmark that is widely used by researchers in the computational geodynamics
mantle convection community; the second demonstrates a global model of the evolution of the Earth's mantle constrained by known movements of the tectonic plates at the surface.

\subsection{Entrainment of a dense layer}
A common benchmark for thermo-chemical mantle convection codes simulates
the entrainment of a dense bottom layer in a 
slow viscous convection driven by heating from below and cooling from above
\cite{JGRB:JGRB11097,TK03,Hansen2000401,Kellogg1881,Tackley1998}. The
benchmark challenges advection schemes because it involves the
tracking of an interface over timespans of
many thousand time steps; field-based methods therefore
accumulate significant numerical diffusion. Results
for this problem can be found in \cite{JGRB:JGRB11097,TK03}.
In our example, we will use a minor modification of
the ``thick layer'' test of \cite{TK03} for 3D models,
but with far higher resolutions and far more particles.  

The model computes finite-element based fluid pressures and velocities in a unit
cube using the incompressible Stokes equations coupled with an
advection-diffusion equation for the temperature. The material is
tracked by assigning a scalar density property to particles, and the
density of the material in a cell is computed as the arithmetic average of all
particles in that cell. $T=1$ and $T=0$ are prescribed as temperature
boundary conditions at the top and bottom. All other sides have
insulating thermal boundary conditions. We use tangential, free-slip
velocity boundary conditions on all boundaries. 

The density property of particles with $z \le 0.4$ is initialized so that their buoyancy ratio
compared to the particles above the interface is $B=-1.0$, where
$B=\delta \rho / (\rho_0 \alpha_0 \delta T)$. 
Here, $\delta \rho$ is the density difference between the two phases,
$\rho_0$ is a reference density of the background material, $\delta T$ is the temperature difference across the domain,
and $\alpha_0$ is the thermal expansivity of the background material.
The other material parameters and gravity are chosen
to result in a Rayleigh number of $10^6$, and -- deviating from
\cite{TK03} -- we choose the initial temperature as
$T_0(x,y,z) = \frac 12 + \frac{1}{10} \sin(\pi z)\cos(\pi x)\cos(\pi
y)$.

Previous studies only showed results up to
resolutions of $64^3$ cells (first order
finite element, and finite volume scheme respectively), with 5-40
particles per cell. Here we reproduce this case, although using
second-order accurate finite elements for the velocity and temperature
(case A), and extend the model to a mesh resolution of $128^3$,
with 40 (case B) and 500 (case C) particles per cell. The last of
these cases results in around 90 million degrees of freedom and
more than a billion particles. The models require approximately 15,000 (A) and
32,000 (B and C) time steps respectively. 

\begin{figure}[tb]
\begin{center}
\includegraphics[width=0.7\textwidth]{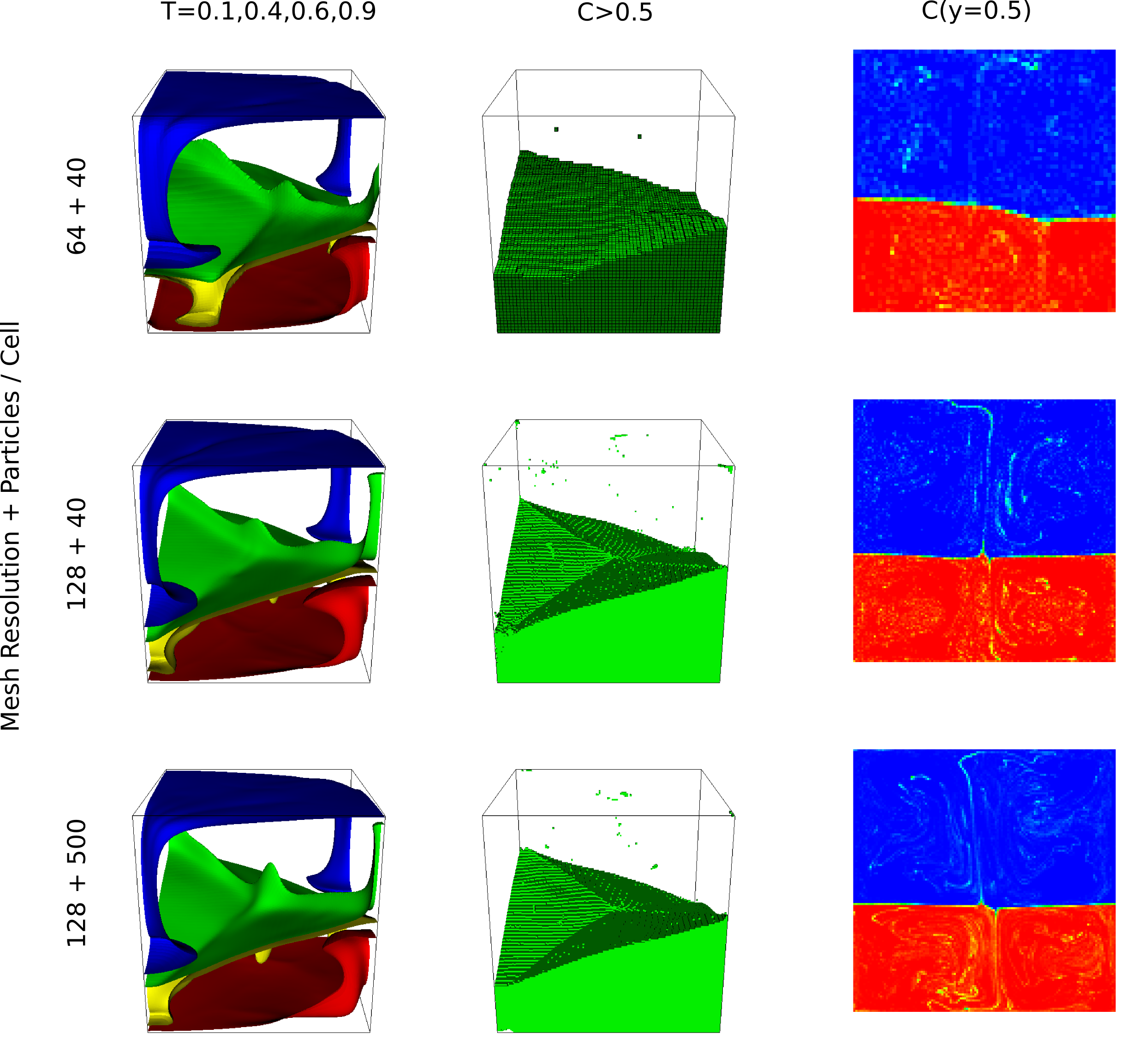}
\end{center}
\caption{\it End-time snapshots of the thick layer test with increasing mesh resolution and particle numbers. Left: Temperature isosurfaces for $T=0.1,0.4,0.6,0.9$. Center: Volume of dense material. Right: Vertical slice, showing the material composition at $y=0.5$.}
\label{fig:benchmark_application}
\end{figure}

Fig.~\ref{fig:benchmark_application} shows the visual end state of these
models, extending Fig.~9 of \cite{TK03}. As in previously
published results, the dense
layer remains stable at the bottom, with small amounts of material
entrained in the convection that develops above it.

Previous studies have measured the entrainment
$
e = \frac{1}{0.4}\int_{0.5}^1 \int_0^1 \int_0^1 C \; \text{d}x\, \text{d}y\, \text{d}z,
$
where $C=1$ for material that originated in the bottom layer, and
$C=0$ otherwise. For case A, we obtain $e=0.02382$, slightly lower than the previously published 3D
result with the same resolution and number of particles ($\sim0.03$). We conjecture
that this is due to using second order elements for velocity
and temperature. Increasing the mesh
resolution to $128^3$ decreases the amount of entrainment further to
$0.01889$, in line with the expected trend \cite{TK03}, but also
showing that a higher resolution is necessary to confirm convergence
for this result. Drastically increasing the number of particles per
cell (case C) only increases the entrainment slightly to $0.01939$,
which is also consistent with the literature. Despite the change
in statistical properties being small, this case resolves more small-scale
features of the flow than case B, in particular the structure of the interface
and the deformation of entrained material.

Thus, these results suggest that the trends for the entrainment rate
reported in
\cite{TK03} for 2D likely also hold for 3D models.
They also show that even though the change in entrainment rate is already low for
few (less than 40) particles per cell, flow structures in the solution are better
resolved by using more than 100 particles per cell.

\subsection{Convection in the Earth's mantle}

Finally, we illustrate the capability of the methodology described above on a 
realistic three-dimensional geodynamics computation.
This problem models the evolution of the Earth's mantle over 250 million years.
The mesh changes adaptively and contains on average 1.7 million cells, resulting in 
90 million degrees of freedom. 4.8 million particles are used for 
post-processing. The particles are initially distributed randomly. 
However, in order to enforce strictly balanced parallel workloads we limit the
maximum number of particles per cell to 25.
Therefore, at the final time those regions of no interest (i.e., those regions resolved with only coarse cells) have a lower particle density than regions of interest.
Specifically, we examine material close to a region of cold downwelling
(a subducted plate or ``slab'') and determine its movement since the beginning of the 
computation.
Material properties such as density and heat capacity are computed from a database for 
basaltic and harzburgitic rocks, following \cite{nakagawa2009}, and the viscosity is 
based on a published viscosity model incorporating mineral physics properties, geoid 
deformation, and seismic tomography \cite{ste06a}. 
The prescribed surface velocities use reconstructions of past plate movement on Earth
\cite{Seton2012}.

\begin{figure}[tb]
\centering
\includegraphics[width=0.7\textwidth]{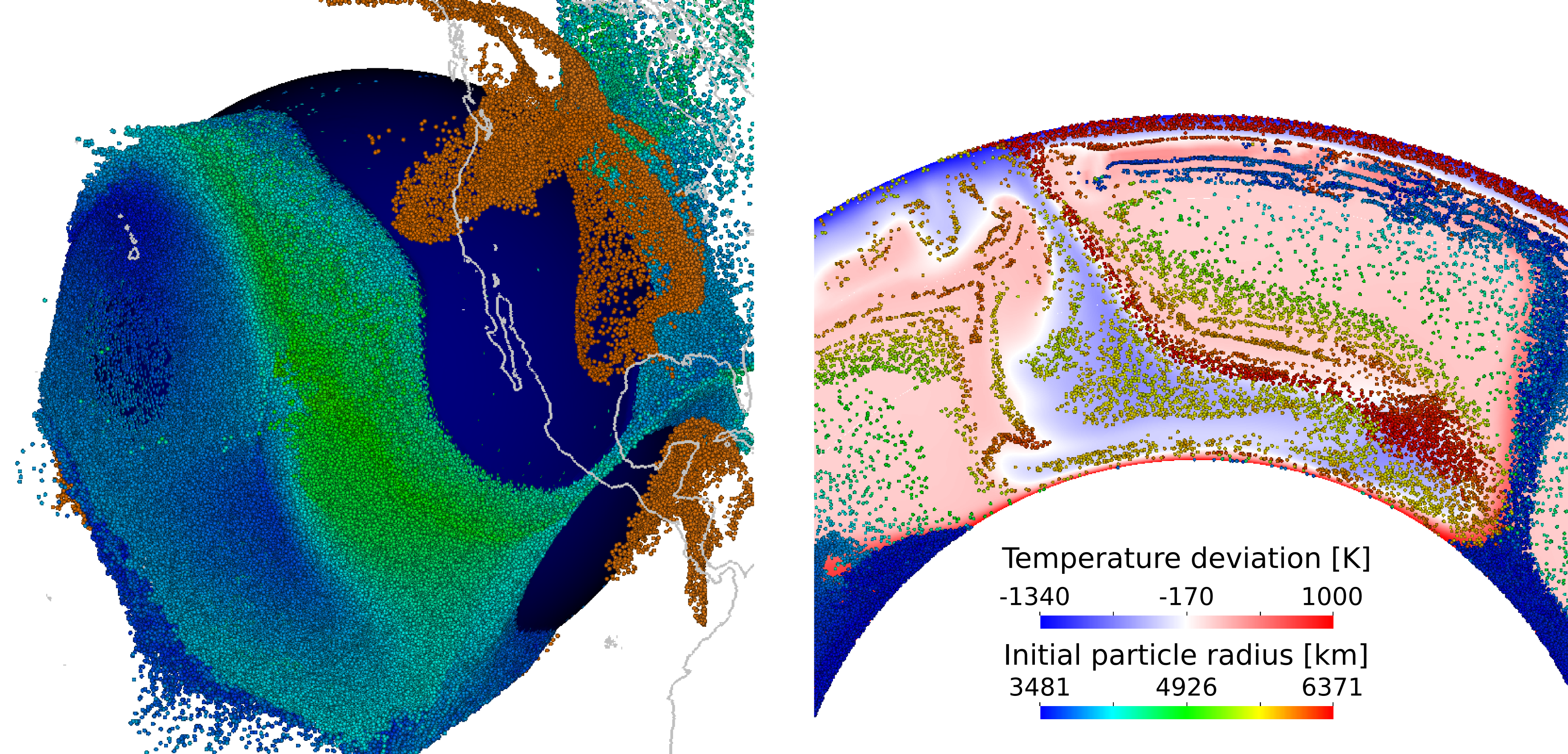}
\caption{\it Illustration of a computation with particles in a 3D mantle
  convection computation. Left: Subducting plates below the western
  United States (orange) push material at the core-mantle boundary
  (dark blue sphere) towards the west. Each Particle is
  colored according to the distance it
  has moved from its initial position (blue: little movement to green: large movement). Right: Vertical slice
  through the subduction zone. The background color is the temperature and
  the particles are colored by their initial distance from the Earth's
  center (red: the particle started at surface; blue: the particle
  started at the core-mantle boundary).}
\label{fig:application}
\end{figure}

Fig.~\ref{fig:application} shows the present-day state of the
Farallon subduction zone below the western United States. 
Particles that are initially close to the core-mantle boundary are colored by the 
displacement they have experienced since the initial time.
This reveals that the Farallon slab (orange) has primarily pushed the easternmost
material.
Particles in the central Pacific have not moved significantly, illuminating the limited
influence of the west Pacific subduction zones.
Studies such as this one enable researchers in computational geodynamics to understand 
the history and dynamics of the Earth's interior.

\section{Conclusions}
\label{sec:conclusion}

In this article, we have reviewed strategies for implementing methods that
couple field- or mesh-based, and particle-based approaches to computational problems in 
continuum mechanics such as fluid flow.
These methods have a long history of use in computational science and engineering, 
However, most approaches for this coupled methodology have only been 
implemented as a sequential code, or as a parallel code in which the domain is 
statically partitioned.
In contrast, modern finite element codes have adaptively refined meshes that change in 
time, have hanging nodes, are dynamically repartitioned as
appropriate, and use complicated communication patterns. 
Consequently, the development of efficient methods to couple continuum and particle
approaches requires a rethinking of the algorithms, the implementation of these 
algorithms, and a consequent update of the available toolset.
We have described how all of the operations one commonly encounters when using particle
methods can be efficiently implemented and have documented through numerical examples 
that the expected optimal complexities can indeed be realized in practice.

\bibliographystyle{plain}
\bibliography{arxiv1}

%\begin{received}
%Received Month Year; revised Month Year; accepted Month Year
%\end{received}

\end{document}